\documentclass[10pt,onecolumn,journal,papersize ]{article}
\usepackage[dvipdfmx]{graphicx}
\def\qed{\hfill \hbox{${\vcenter{\vbox{
   \hrule height 0.4pt\hbox{\vrule width 0.4pt height 6pt
   \kern5pt\vrule width 0.4pt}\hrule height 0.4pt}}}$}}
\newcommand{\maru}[1]{{\ooalign{\hfil#1\/\hfil\crcr
   \raise.167ex\hbox{\mathhexbox20D}}}}
\newcommand{\bysame}{%
   \leavevmode\hbox to 6em{\hrulefill}\,}

\begin{document}

\centerline{\large\bf On Heegaard splittings of knot exteriors with} 
\vskip 2mm

\centerline{\large\bf tunnel number degenerations}
\vskip 1mm 

\centerline{by} 
\vskip 2mm 

\centerline{\bf Kanji Morimoto} 
\vskip 5mm 

\centerline{Department of IS and Mathematics, Konan University} 
\vskip -1mm 

\centerline{Okamoto 8-9-1, Higashi-Nada, Kobe 658-8501, Japan}
\vskip -1mm 

\centerline{morimoto@konan-u.ac.jp} 
\vskip 5mm 

\begin{abstract}

Let $K_1, K_2$ be two knots with $t(K_1)+t(K_2)>2$ and $t(K_1 \# K_2)=2$. Then, in the present paper, we will show that any genus three Heegaard splittings of $E(K_1 \# K_2)$ is strongly irreducible and that $E(K_1 \# K_2)$ has at most four genus three Heegaard splittings up to homeomorphism. Moreover, we will give a complete classification of those four genus three Heegaard splittings and show unknotting tunnel systems of knots $K_1 \# K_2$ corresponding to those Heegaard splittings. 
\vskip 2mm 

\vskip 2mm 

\end{abstract}

\vskip 10mm

{\bf 1. Introduction} 

Let $K$ be a knot in $S^3$ and $t(K)$ the tunnel number of $K$, where $t(K)$ is the minimal number of arcs properly embedded in the exterior $E(K)$ whose complementary space is homeomorphic to a handlebody. By the definition of the tunnel number, we have $t(K)=g(E(K))-1$, where $g(E(K))$ is the Heegaard genus of $E(K)$. 

\hskip 4mm 
Let $K_1$ and $K_2$ be two knots in $S^3$ and $K_1 \# K_2$ the connected sum of $K_1$ and $K_2$. Then, on the degeneration problem of tunnel numbers, i.e., the problem that if there are knots $K_1$ and $K_2$ with $t(K_1 \# K_2)<t(K_1)+t(K_2)$ or not, our first result is the following: 
\vskip 3mm 

{\bf Theorem 1 ([4]).} \it There are infinitely many pairs of knots $K_1$ and $K_2$ such that $t(K_1)=1$, $t(K_2)=2$ and $t(K_1 \# K_2)=2$. \rm 
\vskip 3mm 

\hskip 4mm
Successively, we have characterized such knots as follows: 
\vskip 3mm 

{\bf Theorem 2 ([5]).} \it $(1)$ If $t(K_1)+t(K_2)>2$ and $t(K_1 \# K_2)=2$, then $t(K_1)+t(K_2)=3$. $(2)$ $t(K_1)=1, t(K_2)=2$ and $t(K_1 \# K_2)=2$ if and only if $K_1$ is a $2$-bridge knot and $K_2$ is a knot with a $2$-string essential free tangle decomposition such that at least one of the two tangles has an unknotted component. \rm 
\vskip 3mm 

\hskip 4mm 
In the present paper, we investigate genus three Heegaard splittings of such knot exteriors $E(K_1 \# K_2)$ and show unknotting tunnel systems of $K_1 \# K_2$ corresponding to those Heegaard splittings. First we will show: 
\vskip 3mm 

{\bf Theorem 3.} \it Let $K$ be a tunnel number two knot in $S^3$. Suppose a genus three Heegaard splitting of $E(K)$ is weakly reducible, then $E(K)$ is obtained from $E(K_1)$ and $E_V(K_2)$ by gluing $\partial E(K_1)$ and $\partial V$, where $K_1$ is a tunnel number one knot in $S^3$ and $K_2$ is a tunnel number one knot in a solid torus $V$. \rm  
\vskip 3mm 

\hskip 4mm 
Then we get: 
\vskip 3mm 

{\bf Corollary 1.} \it  Let $K_1$ and $K_2$ be two knots in Theorem $2(2)$. Then any genus three Heegaard splitting of $E(K_1 \# K_2)$ is strongly irreducible. \rm 
\vskip 3mm 

{\bf Remark 1.} In [3], it has been shown by Moriah that genus three Heegaard splittings of $E(K_1 \#K_2)$ are strongly irreducible for some subfamily of those knots $K_1, K_2$ in Theorem 2(2).
\vskip 3mm 

\hskip 4mm 
Next we have: 
\vskip 3mm 

{\bf Theorem 4.} \it Let $K_1$ and $K_2$ be two knots in Theorem $2(2)$. Then $E(K_1 \# K_2)$ has at most four genus three Heegaard splittings up to homeomorophisms. \rm 
\vskip 3mm 
 
\hskip 4mm
To give a complete classification of those four genus three Heegaard splittings in Theorem 4, we assume : 

\hskip 4mm
$K_1$ is a 2-bridge knot $S(\alpha, \beta)$ (Schubert's notation in [10]).

\hskip 4mm
$K_2$ has a 2-string essential free tangle decomposition such that: 

\hskip 10mm $(S^3, K_2) = (C_1, K_2 \cap C_1) \cup (C_2, K_2 \cap C_2)$ and 

\hskip 10mm $C_1$ contains an unknotted component. 

\hskip 4mm
To state the classification theorem, we put the following cases:

\hskip 4mm
Case 1: $C_2$ contains no unknotted component. 

\hskip 4mm
Case 2: $C_2$ contains an unkontted component. 

\hskip 4mm
Furthermore, we divide Case 2 into the following two sub-cases: 

\hskip 4mm
Case 2a: there is a self-homeomorphism of $(S^3, K_2)$ exchanging 

\hskip 18mm 
the two tangles $(C_1, K_2 \cap C_1)$ and $(C_2, K_2 \cap C_2)$.

\hskip 4mm
Case 2b: there is no self-homeomorphism of $(S^3, K_2)$ exchanging 

\hskip 18mm 
the two tangles $(C_1, K_2 \cap C_1)$ and $(C_2, K_2 \cap C_2)$.

\hskip 4mm
Then we get:
\vskip 3mm 

{\bf Theorem 5.} \it Let $K_1$ and $K_2$ be two knots in Theorem $2(2)$. Then we have the following complete classification of genus three Heegaard splittings of $E(K_1 \# K_2)$ up to homeomorphisms, where $n$ is the number of homeomorphism classes. 
\vskip 2mm 

\hskip 4mm 
Case $1$ \ $\left\{ \begin{array}{c} 
n=1 \quad {\rm if} \ \beta \equiv \pm 1 \ ({\rm mod} \ \alpha) \\ 
n=2 \quad {\rm if} \ \beta \not\equiv \pm 1 \ ({\rm mod} \ \alpha) \end{array} \right.$
\vskip 2mm 

\hskip 4mm 
Case $\rm 2a$ $\left\{ \begin{array}{c} 
n=1 \quad {\rm if} \ \beta \equiv \pm 1 \ ({\rm mod} \ \alpha) \\ 
n=2 \quad {\rm if} \ \beta \not\equiv \pm 1 \ ({\rm mod} \ \alpha) \end{array} \right.$
\vskip 2mm 

\hskip 4mm 
Case $\rm 2b$ $\left\{ \begin{array}{c} 
n=2 \quad {\rm if} \ \beta \equiv \pm 1 \ ({\rm mod} \ \alpha) \\ 
n=4 \quad {\rm if} \ \beta \not\equiv \pm 1 \ ({\rm mod} \ \alpha) \end{array} \right.$ \rm 
\vskip 5mm 

{\bf Remark 2.} The condition $\beta \equiv \pm 1 \ ({\rm mod} \ \alpha)$ is equivalent to that $K_1$ is a torus knot.  
\vskip 3mm 

{\bf Example 1.} In Figure 1, (i) is a 2-string essential free tangle with an unknotted component, and (ii) is a 2-string essential free tangle witn no unknotted component. 

\begin{figure}[htbp]
\centerline{\includegraphics[width=6cm]{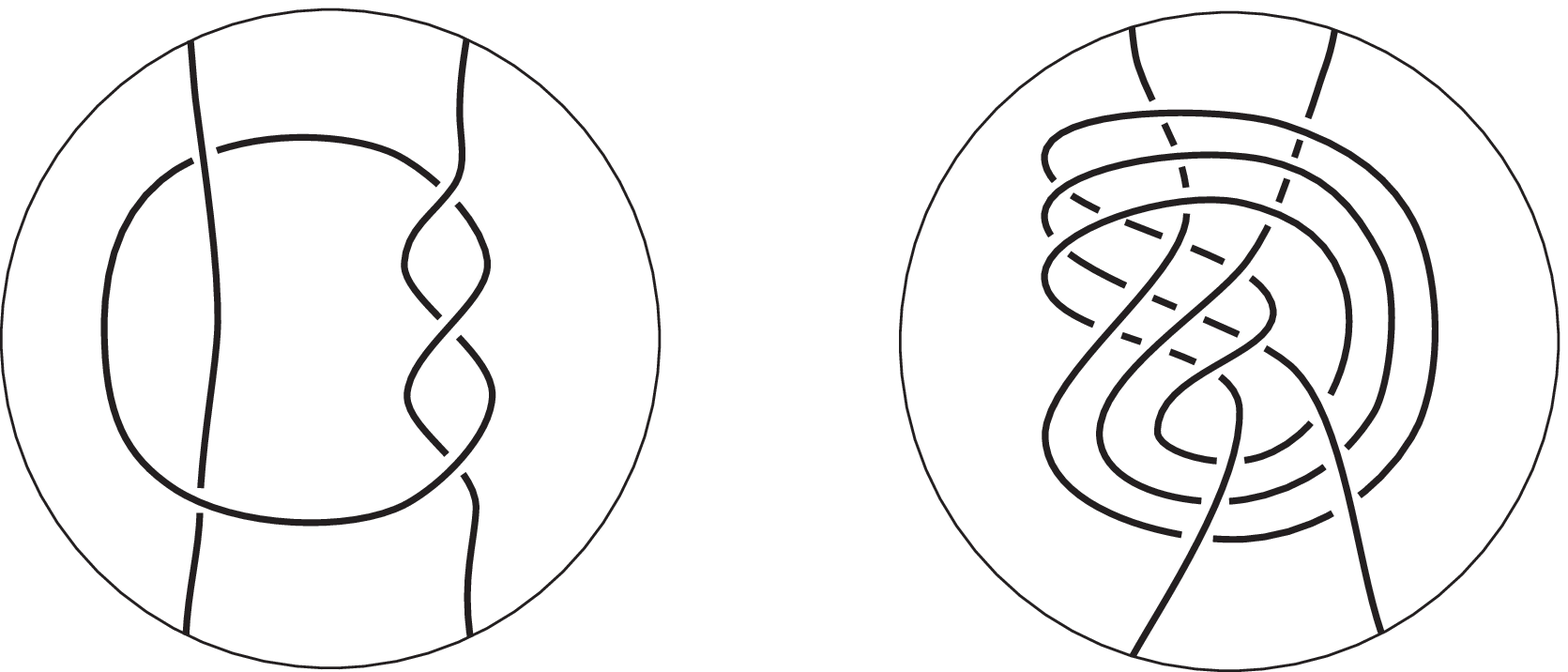}}

\centerline{(i) \hskip 27mm (ii)}

\centerline{Figure 1: 2-string essential free tangles}
\end{figure}

{\bf Example 2.} In Figure 2, (i) is a knot which has a 2-string essential free tangle decomposition such that one of the tangles has an unknotted component, and (ii) is a knot which has a 2-string essential free tangle decomposition such that both tangles have unknotted components, i.e., (i) is in Case 1 and (ii) is in Case 2 of Theorem 5.  

\begin{figure}[htbp]
\centerline{\includegraphics[width=7cm]{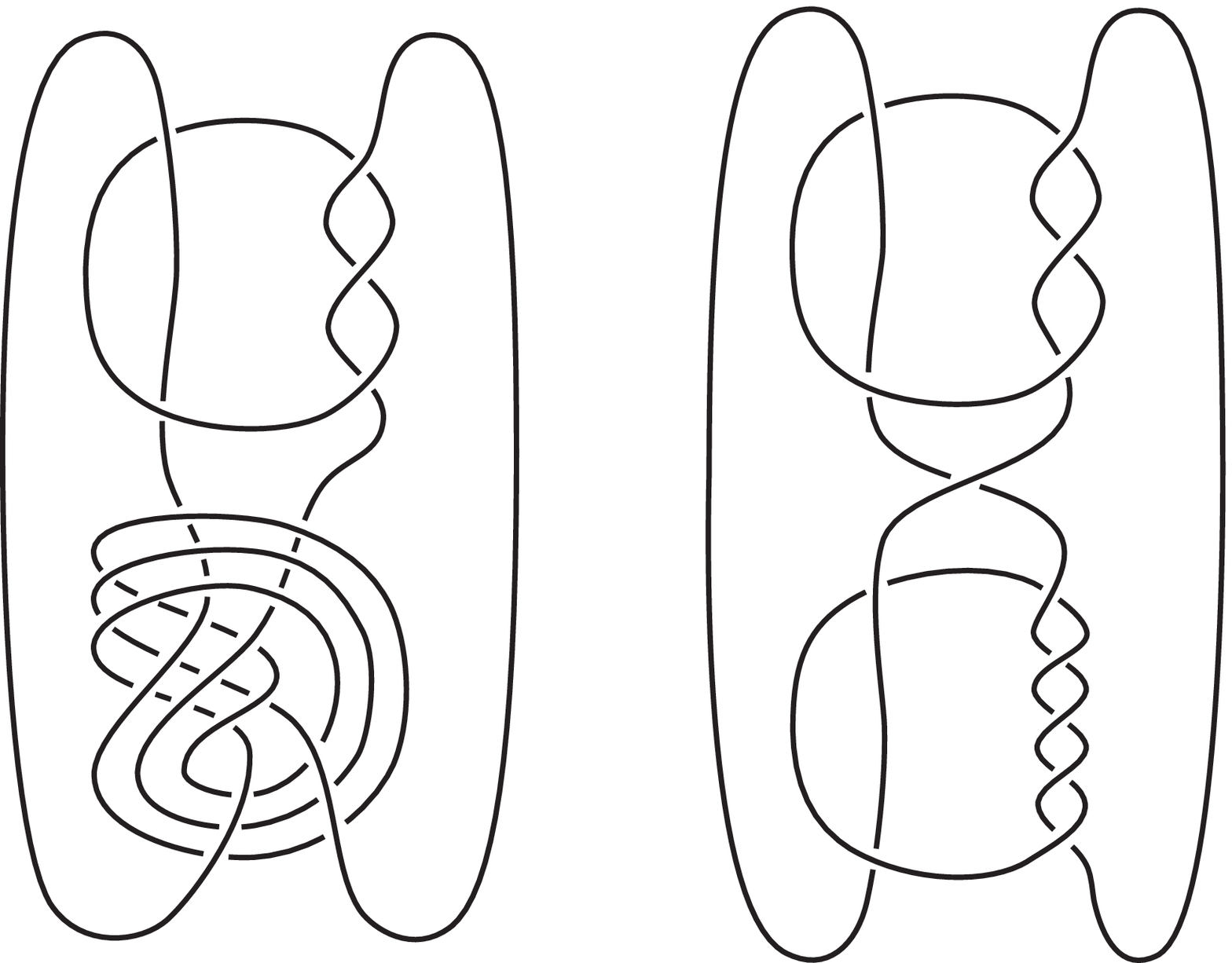}}

\centerline{(i) Case 1 \hskip 20mm (ii) Case 2}

\centerline{Figure 2: Knots with 2-string essential free tangle decompositions}
\end{figure}

{\bf Example 3.} The knot illustrated In Figure 3 is Case 1 of Theorem 5 and the 2-bridge knot is of type (23, 7), i.e., $\beta \not\equiv \pm 1 \ ({\rm mod} \ \alpha)$. Thus the knot exterior of the composite knot has two genus three Heegaard splittings, and the corresponding unknotting tunnel systems are $\{ \tau_1, \tau_2 \}$ and $\{ \sigma_1, \sigma_2 \}$ indicated in the figure. 

\begin{figure}[htbp]
\centerline{\includegraphics[width=5.5cm]{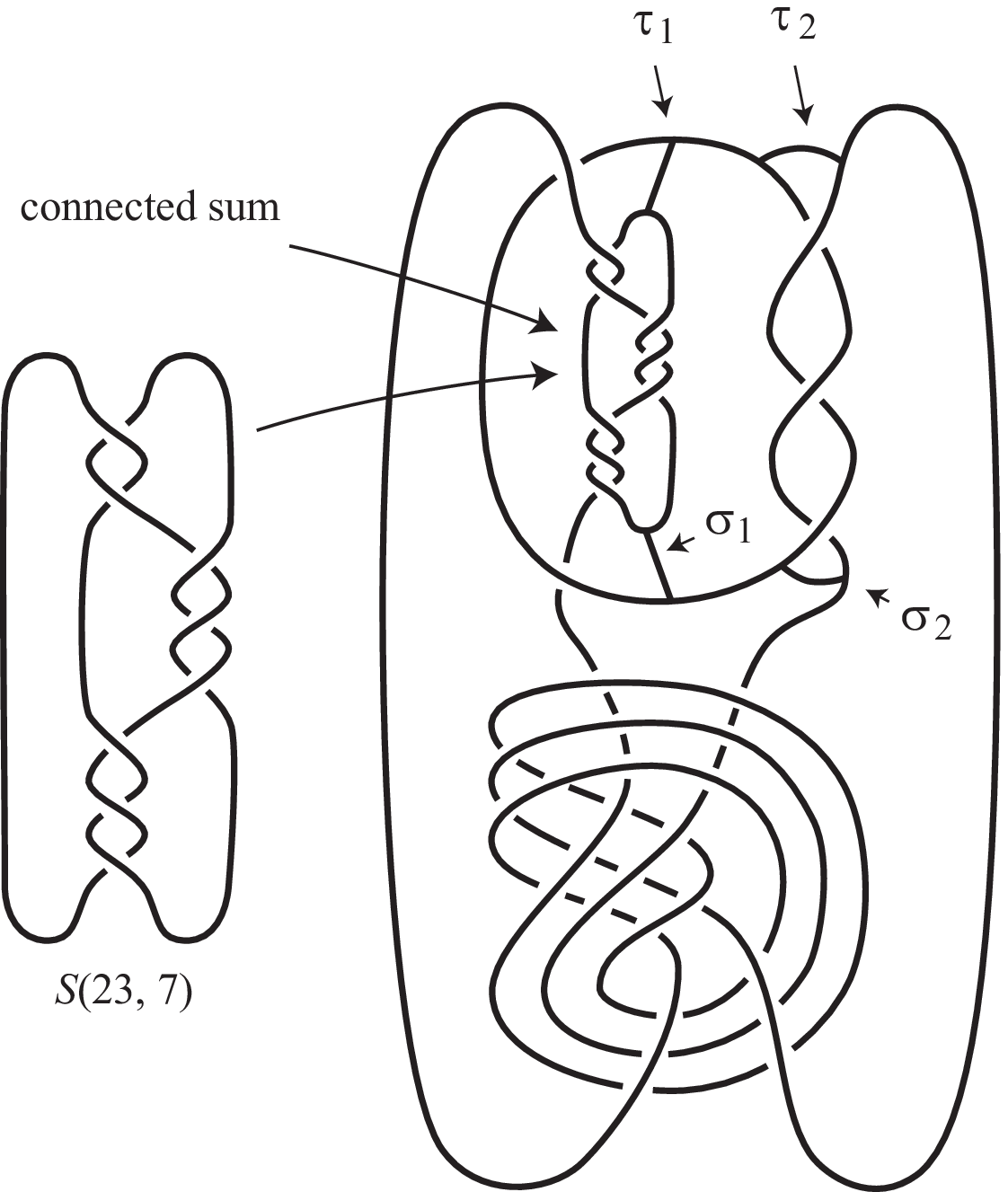}}

\centerline{Figure 3: The two unknotting tunnel systems in Case 1}
\end{figure}

{\bf Example 4.} The two knots illustrated in Figure 4 are the same knots, because by sliding the 2-bridge knot along a sub-arc of the given knot, we can get the right-hand side knot from the left-hand side knot, and this case is Case 2b of Theorem 5. Thus the knot exterior of the knot has four genus three Heegaard splittings and the corresponding unknotting tunnel systems are $\{ \tau_1, \tau_2 \}$, $\{ \sigma_1, \sigma_2 \}$, $\{ \rho_1, \rho_2 \}$ and $\{ \delta_1, \delta_2 \}$ indicated in the figure. 

\begin{figure}[htbp]
\centerline{\includegraphics[width=8cm]{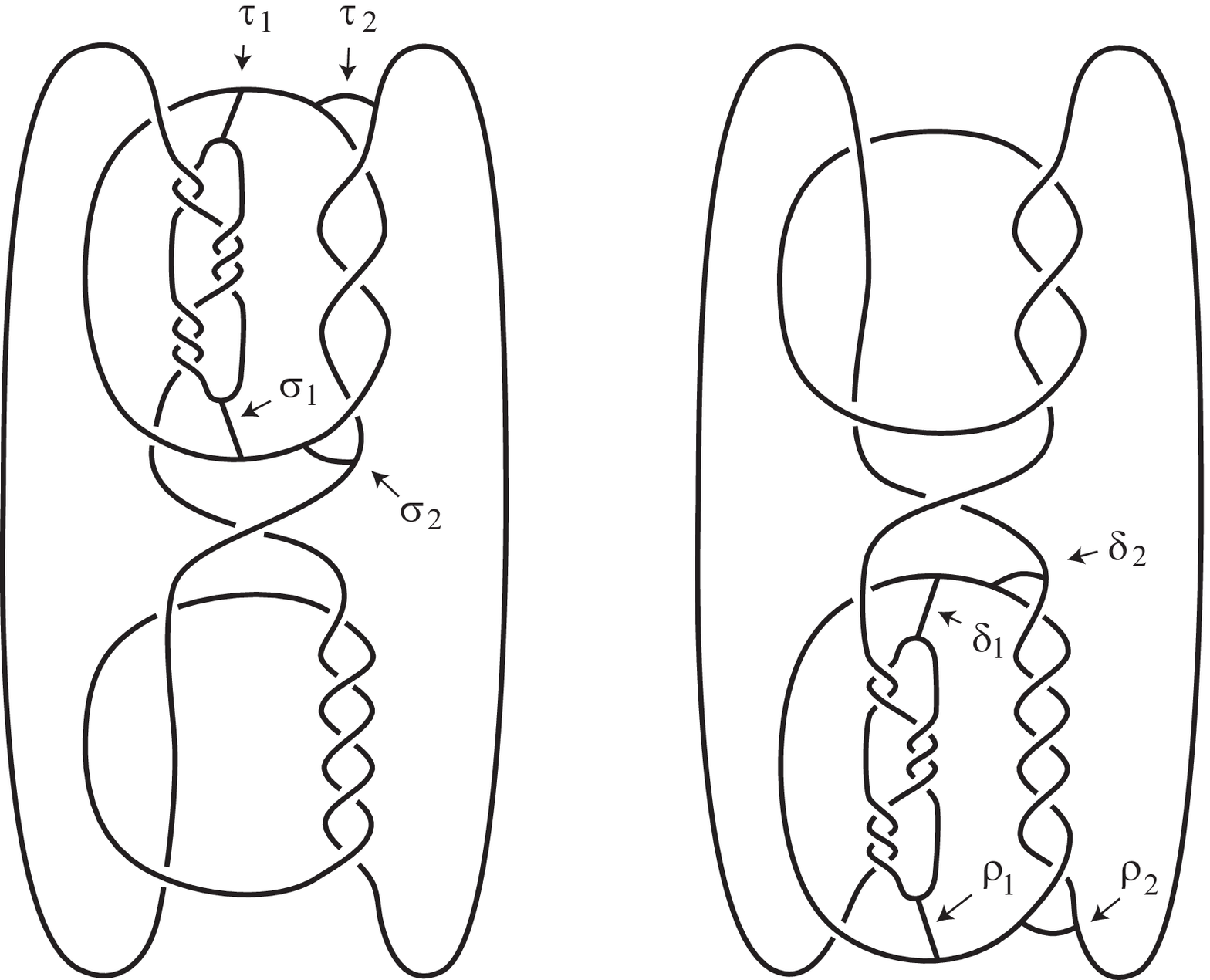}}

\centerline{Figure 4: The four unknotting tunnel systems in Case 2b}
\end{figure}
\vskip 10mm 

{\bf 2. Proofs of Theorem 3 and Corollary 1} 

Let $K$ be a knot in $S^3$, $N(K)$ a regular neighborhood of $K$ in $S^3$ and $E(K)=cl(S^3-N(K))$ the exterior. Put $H_1 \cup H_2$ be a Heegaard splitting of $E(K)$, where $H_1$ is a compression body and $H_2$ is a handlebody, i.e., $\partial E(K) \subset \partial H_1$. We say that the Heegaard splitting $(H_1, H_2)$ is weakly reducible if there is an essential disk, say $D_i$, properly embedde in $H_i$ $(i = 1, 2)$ such that $D_1 \cap D_2 = \emptyset$, and that $(H_1, H_2)$ is strongly irreducible if it is not weakly reducible. For the definition of compression body, we refer [1], and the notion of weak reducibility and strong irreducibility of Heegaard splittings is also due to [1]. 

\hskip 4mm 
Let $V$ be a solid torus and $K$ a knot in $int V$. Let $N_V(K)$ be a regular neighborhood of $K$ in $V$ and $E_V(K) = cl(V - N_V(K))$ the exterior. We say that $K$ is a tunnel number one knot in $V$ if there is an arc $\gamma$ properly embedde in $E_V(K)$ with $\gamma \cap \partial N_V(K) \not= \emptyset$ such that $cl(E_V(K) - N(\partial N_V(K) \cup \gamma))$ is a genus two handlebody (if $\gamma \cap \partial V \not= \emptyset$) or a genus two compression body (if $\partial \gamma \subset \partial N_V(K)$). 
\vskip 3mm 

{\bf Proof of Theorem 3.} Let $H_1 \cup H_2 = E(K)$ be a weakly reducible genus three Heegaard splitting with $\partial E(K) = \partial_{-}H_1$, and $D_1 \subset H_1$ and $D_2 \subset H_2$ be essential disks with $D_1 \cap D_2 = \emptyset$. Then we have the following three cases. 

\hskip 4mm 
Case 1 : Both $D_1$ and $D_2$ are non-separating in $H_1$ and in $H_2$ respectively. 

\hskip 4mm 
Put $H_1' = cl(H_1 - N(D_1))$, $H_2' = cl(H_2 - N(D_2))$, and put $V_1 = cl(H_1' - N(\partial H_1' - \partial E(K)))$, $V_2 = N(\partial H_1' - \partial E(K)) \cup N(D_2)$, $W_1 = N(\partial H_2') \cup N(D_1)$ and $W_2 = cl(H_2' - N(\partial H_2'))$ as illustrated in Figure 5. 

\begin{figure}[htbp]
\centerline{\includegraphics[width=10cm]{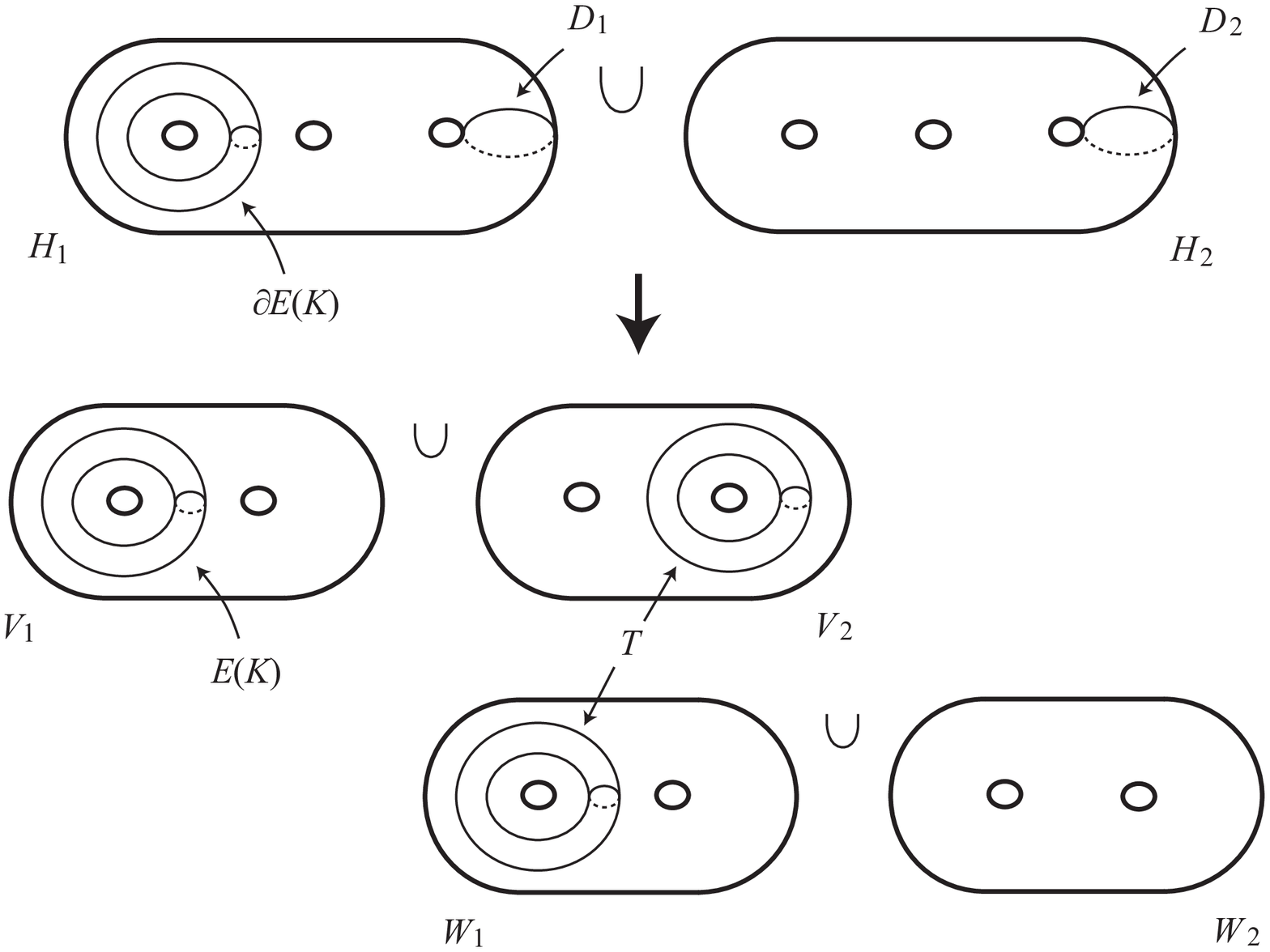}}

\centerline{Figure 5: $D_1$ and $D_2$ are non-separating.}
\end{figure}

\hskip 4mm 
Put $T = V_2 \cap W_1$. Then $T$ is an incompressible torus in $E(K)$ and $(H_1, H_2)$ is an amalgamation of $(V_1, V_2)$ and $(W_1, W_2)$ via $T$. By the solid torus theorem, $T$ is a boundary of a solid torus, say $U$, in the $S^3 = E(K) \cup N(K)$, and $N(K)$ is contained in the sorid torus. Hence $W_1 \cup W_2$ is a knot exterior of some tunnel number one knot in $S^3$ because $(W_1, W_2)$ is a genus two Heegaard splitting. In addition, $V_1 \cup V_2$ is a knot extrior of some tunnel number one knot in the solid torus $U$ because $(V_1, V_2)$ is a genus two Heegaard splitting. 

\hskip 4mm 
Case 2 : Both $D_1$ and $D_2$ are separating in $H_1$ and in $H_2$ respectively. Let $P_i$ be the torus with one hole bounded by $\partial D_i$ in $\partial H_i \ (i = 1,2)$. If $P_1 \cap P_2 \ne \emptyset$, then since $\partial D_1 \cap \partial D_2 = \emptyset$, we have $P_1 \subset P_2$ or $P_2 \subset P_1$. Then by some isotopy, we may assume that $P_1 = P_2$ and $\partial D_1 = \partial D_2$. Then $D_1 \cup D_2$ is a 2-sphere which bounds a 3-ball in $E(K)$. Then the knot $K$ is a trivial knot or a tunnel number one knot, and this is a contradiction. 

\hskip 4mm 
Hence $P_1 \cap P_2 = \emptyset$. Let $T_i = P_i \cup D_i$ be a torus in $H_i \ (i = 1,2)$. If $T_1$ bounds a solid torus in $H_1$, then we can take a meridian disk in the solid torus, and we can take a meridian disk in the solid torus bounded by $T_2$ in $H_2$. Then this case is reduced to Case 1. 

\hskip 4mm 
Suppose $T_1$ bounds a torus $\times I$ in $H_1$, say $X$, and $T_2$ bounds a solid torus in $H_2$, say $Y$. Put $H_1' = cl(H_1 - X)$, $H_2' = cl(H_2 - Y)$, and put $V_1 = cl(H_2' - N(\partial H_2'))$, $V_2 = N(\partial H_2') \cup X$, $W_1 = N(\partial H_1') \cup Y$ and $W_2 = cl(H_1' - N(\partial H_1'))$ as illustrated in Figure 6. 

\begin{figure}[htbp]
\centerline{\includegraphics[width=10cm]{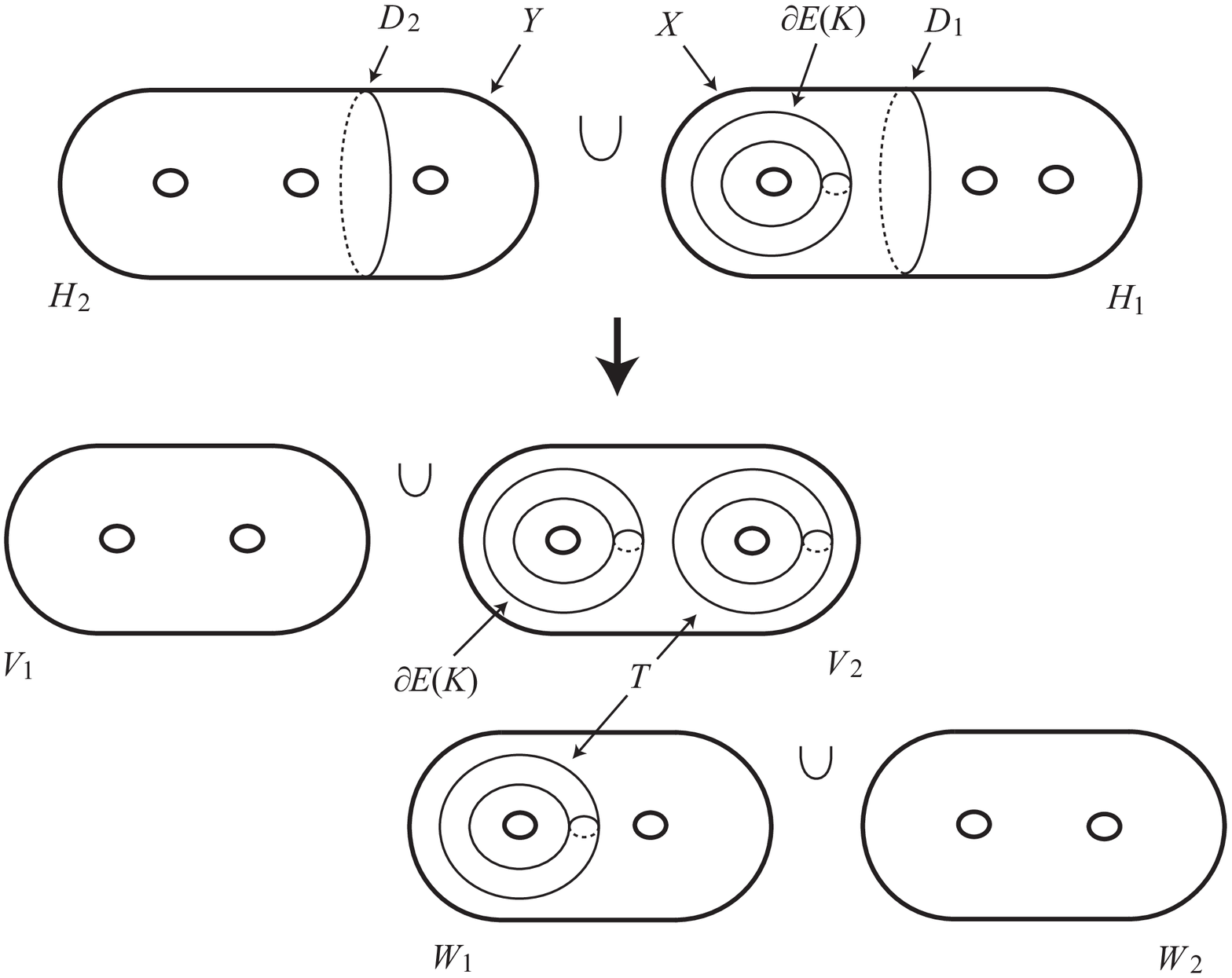}}

\centerline{Figure 6: $D_1$ and $D_2$ are separating.} 
\end{figure}

\hskip 4mm 
Then by the reason similar to the proof of Case I, we see that $W_1 \cup W_2$ is a tunnel number one knot exterior in $S^3$, and $V_1 \cup V_2$ is a tunnel number one knot exterior in a solid torus. 

\hskip 4mm 
Case 3 : One of $D_1$ and $D_2$ is separating and the other is non-separating. 

\hskip 4mm 
Suppose $D_1$ is separating in $H_1$ and $D_2$ is non-separating in $H_2$. Since $\partial D_1 \cap \partial D_2 = \emptyset$, we can take a loop $\ell$ in $\partial H_1 = \partial H_2$ such that $\ell \cap \partial D_1 = \emptyset$ and $\ell \cap \partial D_2$ is a single point. Take a regular neighborhood of $D_2 \cup \ell$ in $H_2$, then it is a solid torus in $H_2$ and let $D_2'$ be the frontier of the solid torus in $H_2$. Then $D_2'$ is a separating essential disk in $H_2$ with $\partial D_1 \cap \partial D_2' = \emptyset$. Next suppose $D_1$ is non-separating in $H_1$ and $D_2$ is separating in $H_2$. Then similarly as above, we can take a separating disk $D_1'$ in $H_1$ with $\partial D_1' \cap \partial D_2 = \emptyset$. Hence Case 3 is reduced to Case 2, and this completes the proof of Theorem 3. \qed 
\vskip 3mm 

{\bf Proof of Corollary 1.} Put $K = K_1 \# K_2$ and suppose $E(K)$ has a genus three weakly reducible Heegaard splitting. Then by Theorem 3, there is an essential torus, say $T$, in $E(K)$ which divides $E(K)$ into a tunnel number one knot exterior in $S^3$, say $E(K_1')$ and a tunnel number one knot exterior in a solid torus $V$, say $E_V(K_2')$. 

\hskip 4mm 
Suppose $T$ is a swallow follow torus of the connected sum. Then, since $t(K_1) = 1$ and $t(K_2) = 2$, $E(K_1)$ is homeomorphic to $E(K_1')$ and $E(K_2)$ is homeomorphic to $E_V(K_2') \cup V'$ for some solid torus $V'$. This shows that $E(K_2)$ has a genus two Heegaard splitting and $t(K_2) = 1$. This is a contradiction., and $T$ is not a swallow follow torus. 

\hskip 4mm 
Let $A$ be the decomposing annulus properly embedded in $E(K)$ corresponding to the connected sum of $K$. 

\hskip 4mm 
First suppose $T \cap A = \emptyset$.  

\hskip 4mm 
If  $T \subset E(K_1)$, then Since $T$ is not a swallow follow torus, $T$ is an essential torus in $E(K_1)$. But 2-bridge knot exterior contains no essential torus by [11]. This is a contradiction.  
If $T \subset E(K_2)$, then by the same reason as above, $T$ is an essential torus in $E(K_2)$. But by [8, Theorem 1.2 and Lemma 1.3] or by [6, Proposition 2.1], this is a contradiction. 

\hskip 4mm 
Hence $T \cap A \ne \emptyset$. Then, since we may assume that each component of $T \cap A$ is an essential loop in both $T$ and $A$, we can take an essential annulus properly embedded in the 2-bridge knot exterior $E(K_1)$ whose boundary components are meridian loops. But this is a contradiction because 2-bridge knots are prime. After all, these contradictions show that $E(K)$ has no genus three weakly reducible Heegaard splitting, and this completes the proof of Corollary 1. \qed 
\vskip 10mm

{\bf 3. Proof of Theorem 4}

Put $K = K_1 \# K_2$, and let $H_1 \cup H_2 = S^3$ be a genus three Heegaard splitting such that $H_1$ contains a knot $K$ as a central loop of a handle of $H_1$. Let $S$ be a decomposing 2-sphere of the connected sum $K_1 \# K_2$. Then by [5], we may assume that $S \cap H_1$ consists of two non-separating disks, say $D_1$ and $D_2$, intersecting $K$ in a single point and a non-separating annulus, say $A$, and that $S \cap H_2$ consists of two non-separating annuli, say $A_1 \cup A_2$, as illustrated in Figure 7. 

\begin{figure}[htbp]
\centerline{\includegraphics[width=12cm]{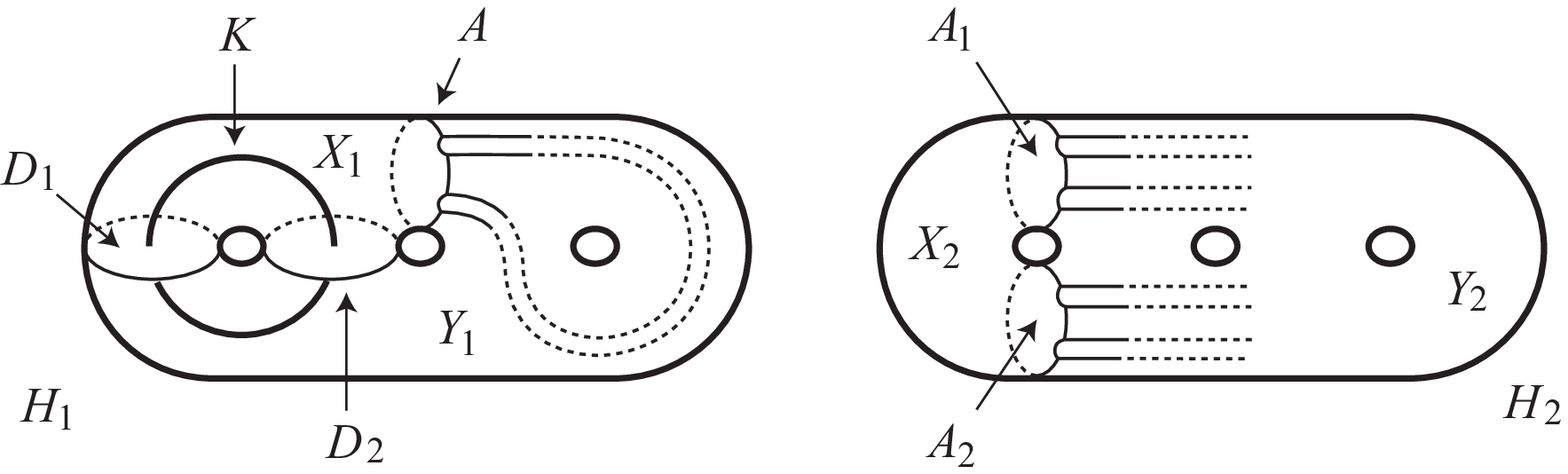}}

\centerline{Figure 7: Heegaard splitting $(H_1, H_2)$} 
\end{figure}

\hskip 4mm 
Then, $S$ splits $H_1$ into two solid tori $X_1$ and $Y_1$ indicated in Figure 7, and $S$ splits $H_2$ into two genus two handlebodies $X_2$ and $Y_2$ indicated in Figure 7. Put $I_1 = [0, 1], I_2 = [1, 2], I_3 = [2, 3]$ and $I = I_1 \cup I_2 \cup I_3$ be intervals, $D_x$ and $D_y$ be two disiks, and P$_x$ and P$_y$ be the central points of $D_x$ and $D_y$ respectively. 

\begin{figure}[htbp]
\centerline{\includegraphics[width=6cm]{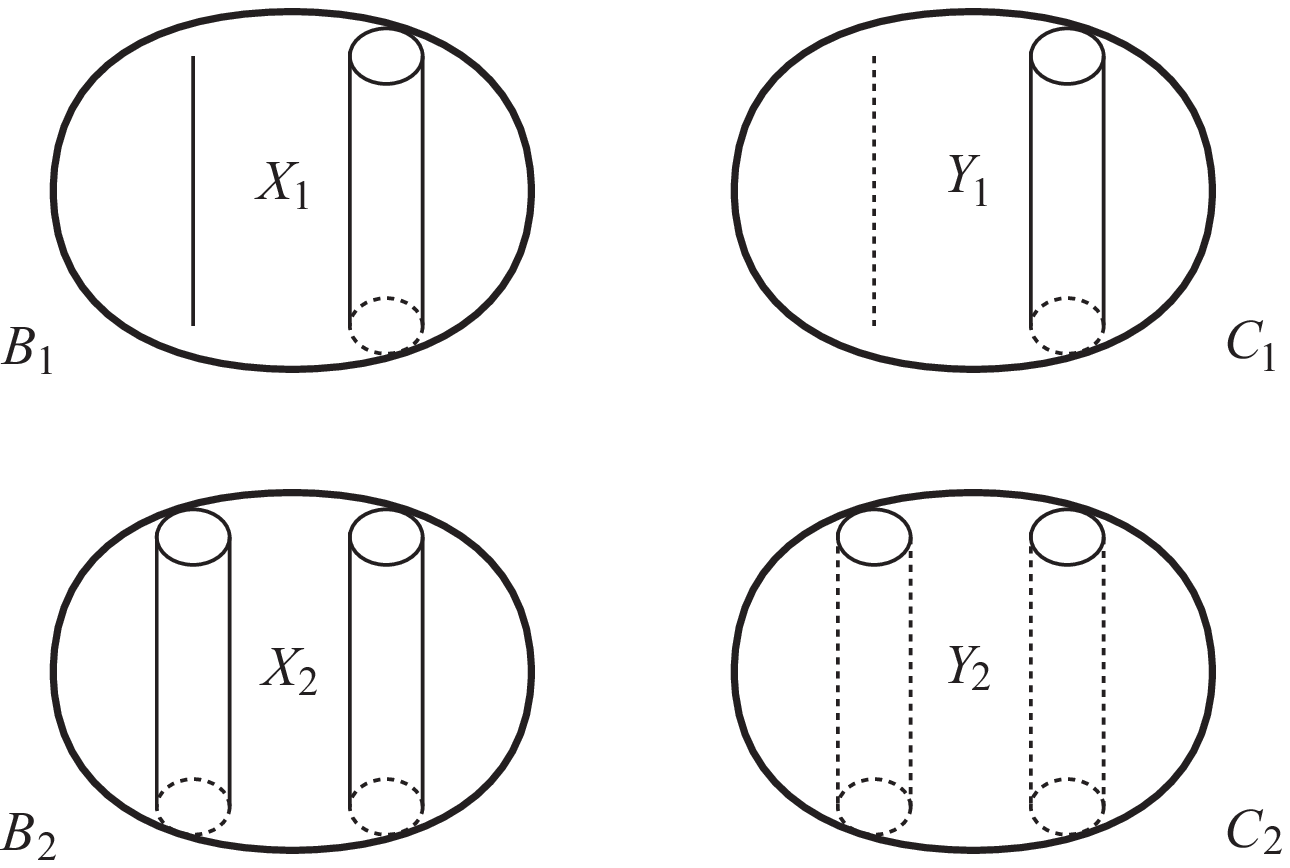}}

\centerline{Figure 8: Tangle decompositions $(B_1, B_2)$ and $(C_1, C_2)$ } 
\end{figure}

\hskip 4mm 
Put $B_1 = X_1 \cup_{(A=\partial D_x \times I_2)} (D_x \times I_2)$ and $B_2 = X_2 \cup_{(A_1 \cup A_2) = (\partial D_x \times (I_1 \cup I_3))} (D_x \times (I_1 \cup I_3))$. Then, since $A$ and $A_1 \cup A_2$ are primitive annuli in $\partial X_1$ and in $\partial X_2$ respectively, $B_1$ and $B_2$ are two 3-balls and $(B_1, B_2)$ gives a 2-bridge decomposition of the knot $K_1 = (B_1 \cap K) \cup ({\rm P}_x \times I)$ in the 3-sphere $B_1 \cup B_2$ (Figure 8). On the other hand, put $C_1 = Y_1 \cup_{(A=\partial D_y \times I_2)} (D_y \times I_2)$ and $C_2 = Y_2 \cup_{(A_1 \cup A_2) = (\partial D_y \times (I_1 \cup I_3))} (D_y \times (I_1 \cup I_3))$. Then, the arguments in the proof of the main theorem of [5] show that both $C_1$ and $C_2$ are 3-balls, and $(C_1, C_2)$ gives a 2-string essential free tangle decomposition of the knot $K_2 = (C_1 \cap K) \cup ({\rm P}_y \times I)$ in the 3-sphere $C_1 \cup C_2$. We note that P$_x \times I_2$ is an unknotted component in $C_1$ (Figure 8).  

\hskip 4mm 
By the above arguments, we can see that any genus three Heegaard splitting of $E(K)$ is obtained from a 2-bridge decomposition of $K_1$ and a 2-string essential free tangle decomposition of $K_2$ by gluing $\partial (D_x \times I) = \partial (X_1 \cup X_2)$ and $\partial (D_y \times I) = \partial (Y_1 \cup Y_2)$. Then, by the uniqueness of prime decomposition of knots ([9]), by the uniqueness of 2-bridge decomposition ([10]), and by the uniqueness of 2-string essential free tangle decomposition ([8]), we have at most four choices of the gluing map up to homeomorphism, i.e., exchanging of $B_1$ and $B_2$ and exchanging of $C_1$ and $C_2$. See (i) $\sim$ (iv) of Figure 9. We note that $X_1', X_2', Y_1'$ and $Y_2'$ in Figure 9 are other components of Heegaard splittings of $E(K)$ (c.f. Figure 10).  Then, by $2 \times 2 = 4$, we complete the proof of Theorem 4. \qed 

\begin{figure}[htbp]
\centerline{\includegraphics[width=12cm]{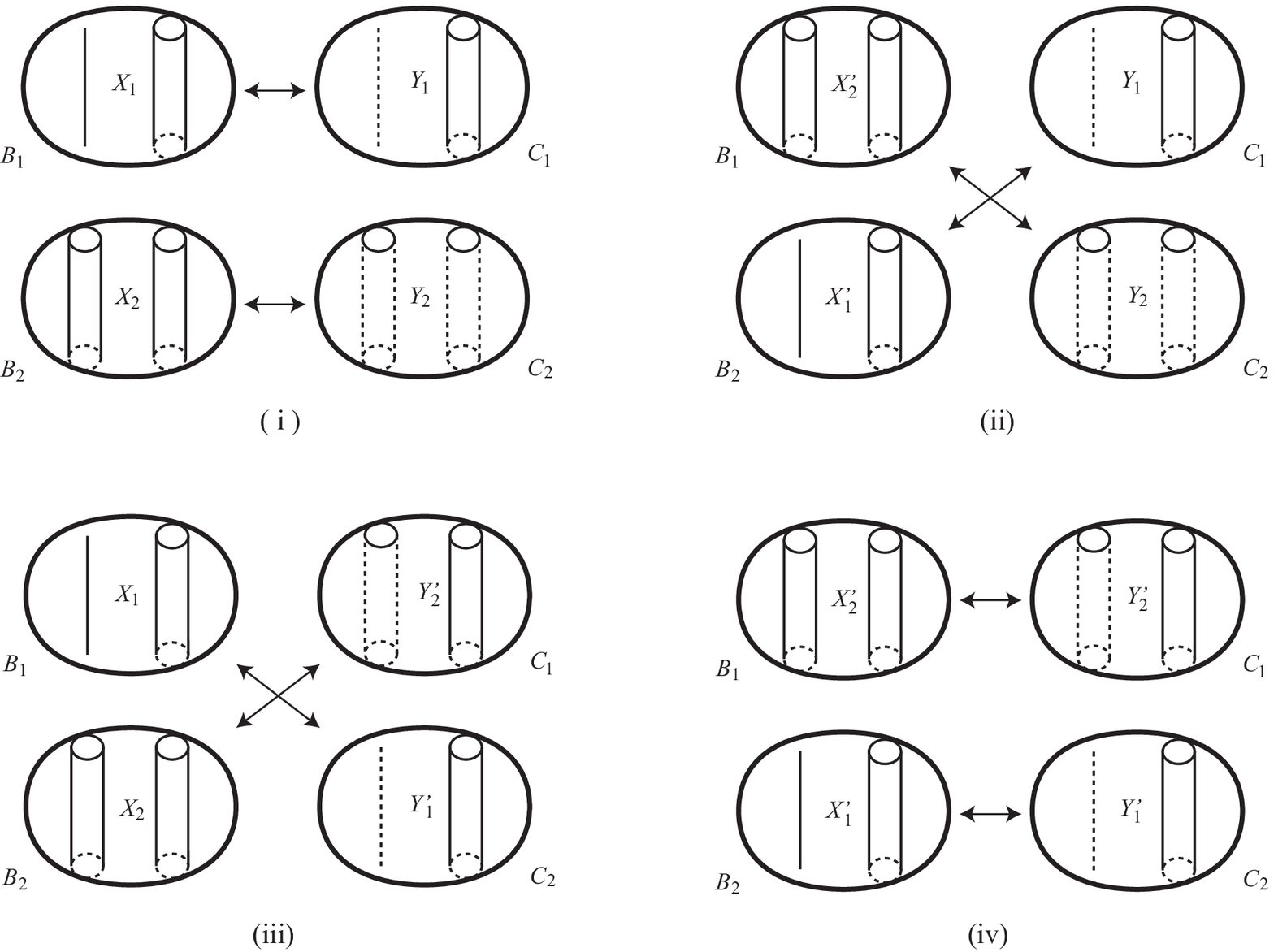}}

\centerline{Figure 9: Four combinations} 
\end{figure}
\vskip 10mm

{\bf 4. Proof of Theorem 5}

As we see the proof of Theorem 4, genus three Heegaard splittings of $E(K)$ are dependent on the choice of 2-bridge decomposition of $K_1$ and free tangle decomposition of $K_2$. 

\hskip 4mm
Suppose we are in Case 1. Then, since $C_2$ contains no unknotted component, we have two Heegaard splittings $(H_1, H_2)$ and $(H_1', H_2')$ such that $H_1 = X_1 \cup Y_1$, $H_2 = X_2 \cup Y_2$, $H_1' = X_1' \cup Y_1$, $H_2' = X_2' \cup Y_2$, where $(X_1, X_2)$ corresponds to $(B_1, B_2)$, $(X_1', X_2')$ corresponds to $(B_2, B_1)$ and $(Y_1, Y_2)$ corresponds to $(C_1, C_2)$. See (i) and (ii) of Figure 9 and Figure 10.  

\begin{figure}[htbp]
\centerline{\includegraphics[width=12.5cm]{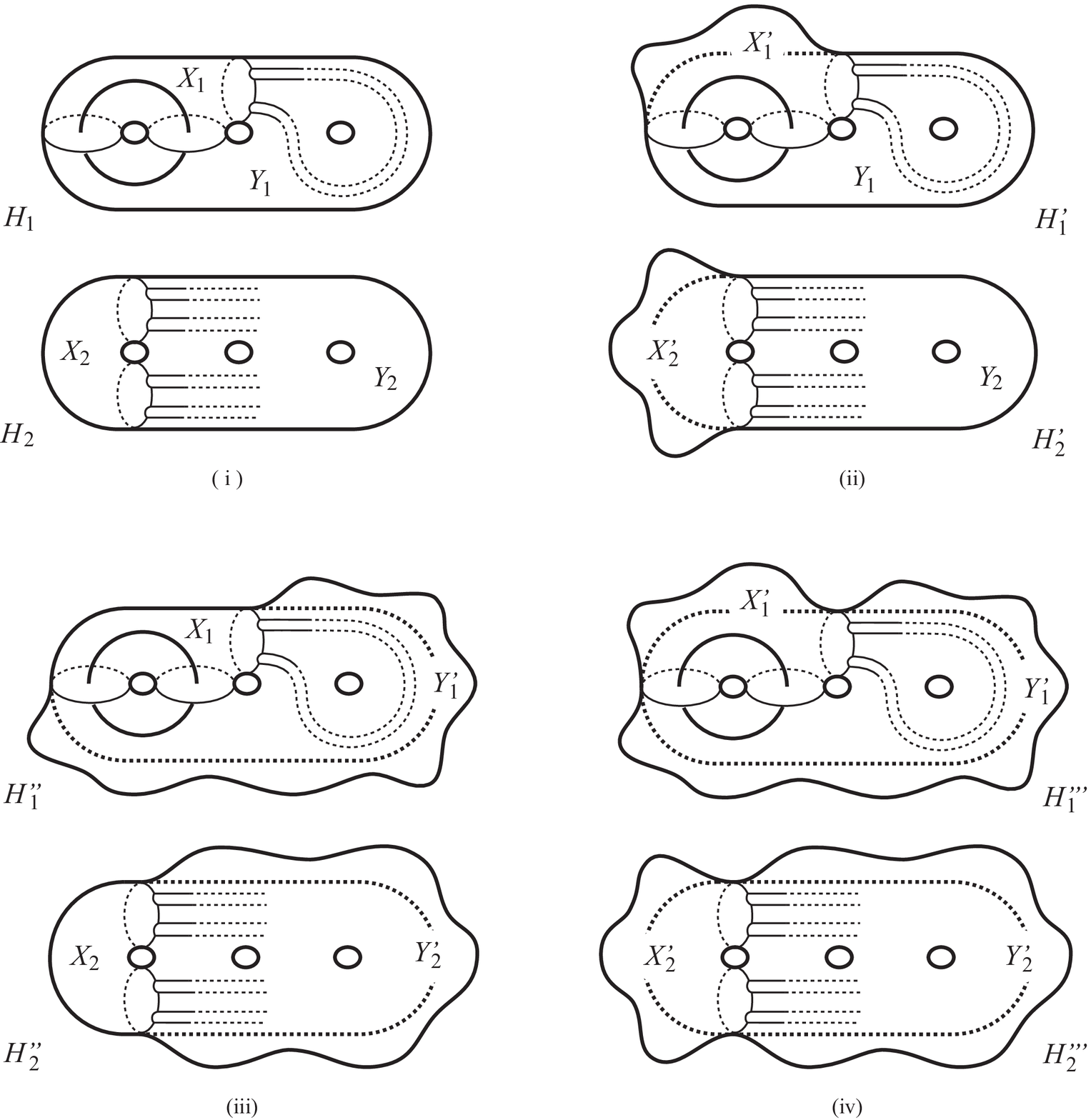}}

\centerline{Figure 10: Four Heegaard splittings} 
\end{figure}

\hskip 4mm 
If $\beta \equiv \pm 1 \ ({\rm mod} \ \alpha)$, then by [7, Theorem 5.2], $B_1$ is isotopic to $B_2$ rel. to $K_1$ in the 3-sphere $S^3 = B_1 \cup B_2$. This implies that $X_1'$ is isotopic to $X_1$ and $X_2'$ is isotopic to $X_2$. Thus $(H_1, H_2)$ is isotopic to $(H_1', H_2')$ and we have $n = 1$. 

\hskip 4mm 
Suppose $\beta \not\equiv \pm 1 \ ({\rm mod} \ \alpha)$, and suppose $(H_1', H_2')$ is homeomorphic to $(H_1, H_2)$. Then the homeomorphism between  $(H_1', H_2')$ and $(H_1, H_2)$ takes $(X_1', X_2')$ to $(X_1, X_2)$, and this homeomorphism induces a self-homeomorphism on  $S^3 = B_1 \cup B_2$ exchanging $B_1$ and $B_2$ rel.to $K_1$. Then, since $\beta \not\equiv \pm 1 \ ({\rm mod} \ \alpha)$ and by [7, Theorem 5.2], this homeomorphism reverses the orientation of the 2-bridge knot $K_1$, and this shows that the homeomorphism between  $(H_1', H_2')$ and $(H_1, H_2)$ exchanges $A_1$ and $A_2$. This means that there is a self-homeomorphism of $Y_2$ which exchanges $A_1$ and $A_2$.  

\hskip 4mm 
Let $a_1$ and $a_2$ be the central loops of $A_1$ and $A_2$ respectively. Then we can regard $(Y_2, a_1 \cup a_2)$ is a genus two Heegaard diagram of $S^3$ because $Y_2$ is a complementary space of a 2-string free tangle and  $a_1$ and $a_2$ are the central  loops of the 2-handles. Then, by taking complete meridian disk system of the genus two handlebody $Y_2$, we have $\pi_1(Y_2) \cong < x, y \ | \ - >$, where $x$ and $y$ correspond to those meridian disks. Then by $a_1$ and $a_2$, we have words $w_1$ and $w_2$ in the letters $x$ and $y$, and we have $\pi_1(S^3) \cong < x, y \ | \ w_1,  w_2 >$. Then, by [2], the representation of $\pi_1(S^3)$ can be deformed into a standard one by a sequence of mutual substitutions. However, this is impossible because $w_1$ and $w_2$ have the same lengths by the existence of a self-homeomorphism of $Y_2$ exchanging $a_1$ and $a_2$. This contradiction shows that  $(H_1', H_2')$ is not homeomorphic to $(H_1, H_2)$, and shows that $n=2$. 

\hskip 4mm 
Next, suppose we are in Case 2. In this case, since $C_2$ also has an unknotted component, 
We have four Heegaard splittings $(H_1, H_2)$, $(H_1', H_2')$, $(H_1'', H_2'')$ and $(H_1''', H_2''')$ such that $H_1 = X_1 \cup Y_1$, $H_2 = X_2 \cup Y_2$, $H_1' = X_1' \cup Y_1$, $H_2' = X_2' \cup Y_2$, $H_1'' = X_1 \cup Y_1'$, $H_2'' = X_2 \cup Y_2'$ and $H_1''' = X_1' \cup Y_1'$, $H_2''' = X_2' \cup Y_2'$, where $(X_1, X_2)$ corresponds to $(B_1, B_2)$ and $(X_1', X_2')$ corresponds to $(B_2, B_1)$, $(Y_1, Y_2)$ corresponds to $(C_1, C_2)$ and $(Y_1', Y_2')$ corresponds to $(C_2, C_1)$. See (i), (ii), (iii) and (iv) of Figure 9 and Figure 10.

\hskip 4mm 
If we are in Case 2b, then, since there is no homeomorphism exchanging $C_1$ and $C_2$, the situation is similar to Case 1 and we see that (iii) and (iv) are not homeomorphic to (i) or (ii). This shows that $n=2$ if  $\beta \equiv \pm 1 \ ({\rm mod} \ \alpha)$ and $n=4$ if $\beta \not\equiv \pm 1 \ ({\rm mod} \ \alpha)$. 

\hskip 4mm 
Suppose we are in Case 2a. Then, since there is a homeomorphism exchanging $C_1$ and $C_2$, we have a homeomorphism which takes $Y_1' \cup Y_2'$ to $Y_1\cup Y_2$ rel. to $Y_1 \cap K = Y_1' \cap K$ respectively. This homeomorphism induces a self-homeomorphism on $A_1 \cup A_2$ and on $A \cup D_1 \cup D_2$. Then, since any 2-bridge knot is strongly invertible, this homeomorphism extends to a homeomorphism $X_1' \cup X_2'$ to $X_1\cup X_2$ rel. to $X_1 \cap K = X_1' \cap K$ respectively. Thus, this case is reduced to Case 1, and we have $n=1$ if  $\beta \equiv \pm 1 \ ({\rm mod} \ \alpha)$ and $n=2$ if $\beta \not\equiv \pm 1 \ ({\rm mod} \ \alpha)$.  This completes the proof of Theorem 5. \qed 
\vskip 10mm

{\bf 5. Unknotting tunnel systems} 

In the present section, we will show the unknotting tunnel systems corresponding to those Heegaard splittings of Theorem 5. Recall the Heegaard splitting $(H_1, H_2)$ and consider the unknotting tunnel system $\{ \tau_1, \tau_2 \}$ in $H_1$ as in Figure 11. Then $\tau_1$ is divided by $S$ into two arcs $\tau_1' \cup \tau_1''$. Then $\tau_1'$ is an upper or a lower tunnel of the 2-bridge knot $K_1$, $\tau_1''$ is an arc in $C_1$ connecting $K_2 \cap C_1$ and $A$ and $\tau_2$ is a core loop of the solid torus $Y_1$ together with a sub-arc of $K_2$. Then, by applying these situations to the knots $K_1$ and $K_2$ as illustrated in Figure 3 and Figure 4, we have those unknotting tunnel systems illustrated in Figure 3 and Figure 4. In fact, by the deformation (i) $\sim$ (iv) as in Figure 12, we see that $\tau_2$ is in the position of 
Figure 3.    

\begin{figure}[htbp]
\centerline{\includegraphics[width=7cm]{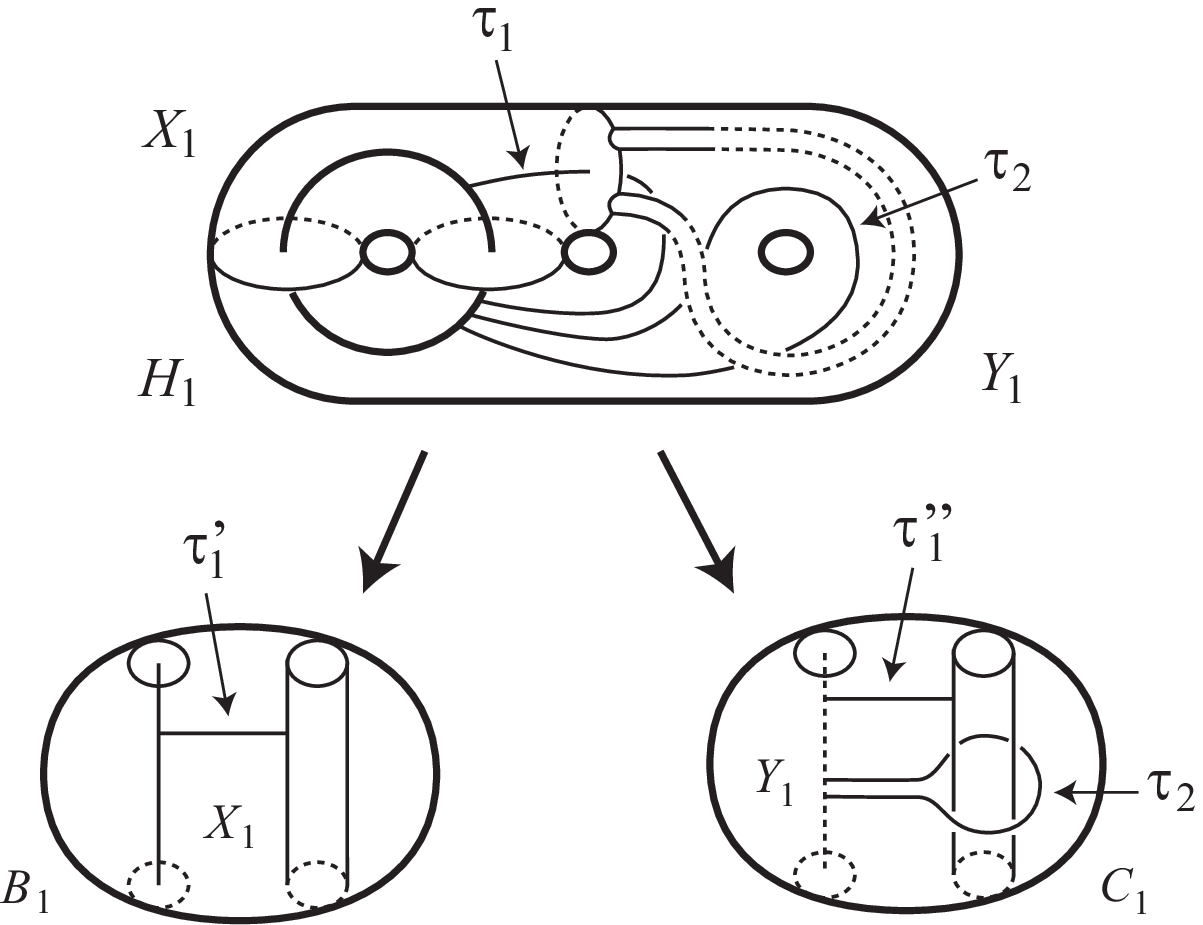}}

\centerline{Figure 11: Heegaard splitting and unknotting tunnel system} 
\end{figure}

\begin{figure}[htbp]
\centerline{\includegraphics[width=11cm]{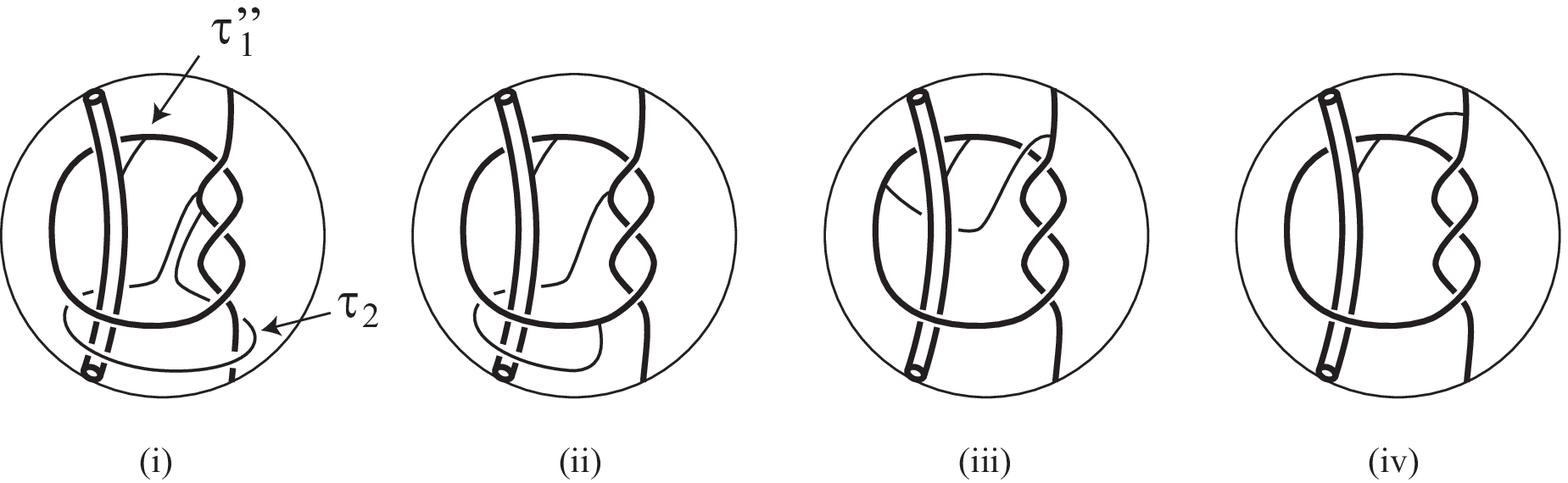}}

\centerline{Figure 12: Deformation of $\tau_2$} 
\end{figure}

\hskip 4mm 
At the end of the present paper, we put the following problem: 
\vskip 3mm 

{\bf Problem} Classify the genus three Heegaard splittings of $E(K_1 \# K_2)$ up to isotopies. 
\vfill\eject

{\bf\large References}
\vskip 3mm 

\ [1] \ A. J. Casson and C. McA. Gordon, Reducing Heegaard splittings, Topology Appl. 

\hskip 7mm {\bf 27} (1987), 275-283.

\ [2] \ T. Kaneto, On genus 2 Heegaard diagram for the 3-sphere, Trans. A. M. S., {\bf 276}

\hskip 7mm (1983), 583-597. 

\ [3] \ Y. Moriah, Connected sums of knots and weakly reducible Heegaard splittings, 

\hskip 7mm Topology Appl. {\bf 141} (2004), 1-20. 

\ [4] \ K. Morimoto, There are knots whose tunnel numbers go down under connected 

\hskip 7mm sum, Proc. A. M. S., {\bf 123} (1995), 3527-3532. 

\ [5] \ K. Morimoto, Charaterization of tunnel number two knots which have the prop- 

\hskip 7mm erty  `` $2 + 1 = 2$ '', Topology Appl. {\bf 64} (1995), 165-176

\ [6] \ K. Morimoto, Essntial tori in 3-string free tangle decmpositions of knots, Journal  

\hskip 7mm of Knot Theory and its Ramifications, {\bf 15}  (2006), 1357-1362. 

\ [7] \ K. Morimoto and M. Sakuma, On unknotting tunnels for knots, Mathematische 

\hskip 7mm Annalen, {\bf 289} (1991), 143-167.

\ [8] \ M. Ozawa, On uniqueness of essential tangle decompositions of knots with free 

\hskip 7mm tangle decompositions, Proc. Appl. Math. Workshop {\bf 8}, ed G.T.Jin and K.H.Ko, 

\hskip 7mm KAIST, Taejon (1998) 

\ [9] \ H. Schubert, Die eindeutige Zerlegbarkeit eines Knoten in Primknoten, Sitzungs- 

\hskip 7mm ber. Akad. Wiss. Hidelberg, math.-nat. KI. {\bf 3} (1949) Abh: 57-104.

[10] \ H. Schubert, Knoten mit zwei Brucken, Math. Z. {\bf 65} (1956), 133-170.

[11] \ H. Schubert, $\ddot U$ber eine numerische Knoteninvariante, Math. Z. {\bf 61} (1954), 245-

\hskip 7mm 288.

\end{document}